\documentclass[12pt,leqno]{amsart}
\usepackage{amssymb}
\usepackage{amsmath,amssymb,color}
\numberwithin{equation}{section}

\newcommand{\weight}{e^{2s\varphi}}

\newcommand{\ep}{\varepsilon}
\newcommand{\la}{\lambda}
\newcommand{\va}{\varphi}
\newcommand{\ppp}{\partial}

\newcommand{\www}{\widehat}

\newcommand{\R}{\mathbb{R}}

\newcommand{\ooo}{\overline}
\newcommand{\OOO}{\Omega}

\newcommand{\sumij}{\sum_{i,j=1}^d}

\newcommand{\hhalf}{\frac{1}{2}}

\allowdisplaybreaks
\title
[]
{
Stability for backward problems in time for degenerate parabolic equations
}

\author{
$^1$ Piermarco Cannarsa and $^2$ Masahiro Yamamoto}

\thanks{
$^1$ Dipartimento di Matematica, Universit\`a di Roma Tor Vergata,\\
Via della Ricerca Scientifica, 00133, Roma, Italy\\
e-mail: {\tt cannarsa@mat.uniroma2.it}\\
$^2$
Graduate School of Mathematical Sciences, The University of Tokyo\\
3-8-1 Komaba, Meguro-ku, Tokyo, 153-8914, Japan.\\
e-mail: {\tt myama@ms.u-tokyo.ac.jp}
}

\date{}

\begin{document}
\maketitle
\begin{abstract}
For solution $u(x,t)$ to degenearte parabolic equations in a bounded domain 
$\OOO$ with homogenous boundary condition, we consider backward problems 
in time: determine $u(\cdot,t_0)$ in $\OOO$ by $u(\cdot,T)$, where $t$ is the 
time variable and $0\le t_0 < T$.
Our main results are conditional stability under boundedness assumptions
on $u(\cdot,0)$.  The proof is based on a weighted $L^2$-estimate of $u$ whose
weight depends only on $t$, which is an inequality of Carleman's type.
Moreover our method is applied to semilinear degenerate parabolic equations.
\end{abstract}

\maketitle

\section{Introduction}

The backward problem for parabolic equations is a typical ill-posed
problem, but it is practically meaningful for applications such as 
estimation of past temperature by means of current temperature distribution, 
which is related to thermo-archaelogy or adequate policies for 
the global warming by estimating the past temperature.

Therefore the mathematical analysis is demanded and the uniqueness and 
the stability are main theoretical topics.
The backward problem in time for the parabolic equation is severelly 
ill-posed, but it is known that we can recover stability under a priori 
boundedness condition.
This recovered stabiliy is called conditional stability and a priori 
boundness recovering the stability, can be introduced from physical viewpoints.
For example, in the heat conduction, such an a priori boundedness 
can be a fusion point of the material under consideration.
Moreover, the parabolic equations of degenerate type are physically
important, and we can refer to enormous references, but here only to
the monograph by  
Cannarsa, Martinez and Vanconstanoble \cite{CMV}, which describes 
physical backgrounds of the degenerate parabolic equations and also 
Carleman estmates with the applications to inverse problems
with comprehensive references.

The main purpose of this article is to establish the conditional stability for 
degenerate parabolic equations.
To the best knowledge of the authors, there are no works on the backward 
problems for degenerate parabolic equations which aim at comprehensive and
systematic researches.

Now we formulate our problem.  Let $\OOO\subset \R^d$ be a bounded domain with 
smooth boundary $\ppp\OOO$.
We consider a degenerate parabolic equation with boundary condition:
$$
(Lu)(x,t):=\ppp_tu(x,t) - \sumij \ppp_i(a_{ij}(x,t)\ppp_ju(x,t)) 
- \sum_{k=1}^d b_k(x,t)\ppp_ku(x,t) 
- c(x,t)u(x,t)
$$
$$
 = F(x,t)         \eqno{(1.1)}
$$
for $(x,t) \in Q:= \OOO\times (0,T)$, and
$$
u=0  \quad \mbox{on $\ppp\OOO \times (0,T)$}    \eqno{(1.2)}
$$
or 
$$
\ppp_{\nu_A}u + r(x)u = 0 \quad \mbox{on $\ppp\OOO
\times (0,T)$}.                                  \eqno{(1.3)}
$$
Here $\nu(x) = (\nu_1(x), ...., \nu_d(x))$ denotes the unit outward normal
vector to $\ppp\OOO$ and we set 
$$
\ppp_{\nu_A}v(x) = \sumij a_{ij}(x)\nu_i\ppp_jv(x) \quad 
\mbox{for $x \in \ppp\OOO$},
$$
and we assume that $r \in C(\ppp\OOO)$.
We assume that 
$$
a_{ij} \in C^1([0,T]; L^{\infty}(\Omega)), \quad
a_{ij} = a_{ji}, \quad 1\le i,j \le d,                  \eqno{(1.4)}
$$
and there exists a function $\sigma \in C(\ooo{\OOO})$, $\ge 0$ on 
$\ooo{\OOO}$ such that 
$$
\sumij a_{ij}(x,t)\xi_i\xi_j \ge \sigma(x)\sum_{k=1}^d \xi_k^2,
\quad (x,t) \in Q, \, \xi_1, ..., \xi_d \in \R,       \eqno{(1.5)}
$$
and $b:= (b_1, ..., b_d) \in L^{\infty}(Q)$,
$c\in L^{\infty}(Q)$.

We set  
$$
H^{2,1}(Q) := \left\{ u\in L^2(Q);\,
\ppp_tu, \ppp_iu, \ppp_i\ppp_ju \in L^2(Q), \quad 1\le i,j\le d
\right\}.
$$
Let $u \in H^{2,1}(Q)$ satisfy (1.1).
Then we consider
\\
{\bf Backward problem.}
\\
{\it 
Let $0\le t_0 < T$ be given.
Determine $u(x,t_0)$, $x\in \OOO$ by 
$u(x,T)$, $x\in \OOO$.
}
\\

The conditional stability has been studied well for the 
non-degenerate parabolic equation where 
$$
\sigma(x) > 0 \quad \mbox{on $\ooo{\OOO}$}.  \eqno{(1.6)}
$$
For the case (1.6), 
as available methodologies, we refer to 
\begin{itemize}
\item
Logarithmic convexity: Ames and Straughan \cite{AS}, Payne \cite{Pay},
Chapter 3 in Isakov \cite{Is}, for example.  
\item
the time analyticity and the maximum principle for holomorphic functions:
Kre{\u\i}n and Prozorovskaya \cite{KP}.
\item
Weight energy methods: \cite{AS}, Lees and Protter \cite{LePr}, 
Payne \cite{Pay}.
\end{itemize}

Our main purpose of this article  
is to establish the conditional stability for the backward problem in time
for the degenerate case 
$$
\sigma(x) \ge 0,\quad x\in \ooo{\OOO},          \eqno{(1.7)}
$$
which means that $\sigma(x)$ in (1.5) is admitted to have zeros.

We introduce the main assumptions on the degeneracy:
$$
\mbox{There exists a constant $\la_1 >0$ such that}\\
$$
$$
\sumij (\la_1a_{ij}(x,t) - \ppp_ta_{ij}(x,t))\xi_i\xi_j
\ge 0, \quad (x,t) \in \ooo{Q}, \, \xi_1, ..., \xi_d \in \R. 
                                                        \eqno{(1.8)}
$$
and
$$
\mbox{there exists a constant $C >0$ such that}\quad
\vert b(x,t)\vert \le C\sqrt{\sigma(x)}, \quad (x,t)\in \ooo{Q}.
                                  \eqno{(1.9)}
$$

If $\vert b(x,t)\vert \equiv 0$ in $Q$, then (1.9) is automatically 
satisfied.

{\bf Examples of (1.8) - (1.9).}
\\
(a) If $a_{ij}(x)$, $1\le i,j \le d$ are 
$t$-independent, then (1.8) is satisfied.
\\
(b) Let $b \equiv 0$ in $Q$.
We assume that we can choose $\mu_k \in C(\ooo{\OOO})$, $\mu_k\ge 0$ on 
$\ooo{\OOO}$, $\www{a_{ij}^k} \in C^1([0,T];L^{\infty}(\OOO))$,
$1\le i,j \le d$, such that there exists a constant $\sigma_1 > 0$ such that 
$$
a_{ij}(x,t) = \sum_{k=1}^N \mu_k(x)\www{a_{ij}^k}(x,t), 
$$
$$
\sumij \www{a_{ij}^k}(x,t)\xi_i\xi_j \ge \sigma_1\sum_{j=1}^d
\xi_j^2 \quad (x,t) \in \ooo{Q}, \, \xi_1, ..., \xi_d \in \R, \quad
1\le k \le N.
                                                          \eqno{(1.10)}
$$
Then (1.8) - (1.9) are satisfied.  For example, 
$\www{a_{ij}^k}(x,t)=\delta_{ij} := 
\left\{ \begin{array}{rl}
1 \quad & i=j,\\
0 \quad & i\ne j
\end{array}\right.$ for $1\le k \le N$ and 
$\mu_k(x) = r_k\vert x-x_0^k\vert^{\rho_k}$ with 
$x_0^k\in \OOO$ and $r_k > 0$, $\rho_k>0$ are constants for 
$1\le k \le N$.
\\
(c) Let $\mu \in C(\ooo{\OOO})$, $\ge 0$ on $\ooo{\OOO}$, 
$\www{a_{ij}} \in C^1([0,T];L^{\infty}(\OOO))$, $b_k\in L^{\infty}(Q)$,
$1\le i,j,k \le d$,
$a_{ij}(x,t) = \mu(x)\www{a_{ij}}(x,t)$ and 
$b_k(x,t) = \mu(x)\www{b_k}(x,t)$ for $1\le i,j,k \le d$.
We further assume that we can find a constant 
$\sigma_1>0$ such that 
$$
\sumij \www{a_{ij}}(x,t) \xi_i\xi_j \ge \sigma_1\sum_{k=1}^d
\xi_k^2 \quad (x,t) \in \ooo{Q}, \, \xi_1, ..., \xi_d \in \R.
$$
Then (1.8) is satisfied.

Now we are ready to state our main results.
\\
\vspace{0.2cm}
{\bf Theorem 1 (case $0<t_0<T$).}
\\
{\it
We assume (1.8) and (1.9), and $u\in H^{2,1}(Q)$ satisfy (1.1) and
$$
\Vert u(\cdot,0)\Vert_{H^1(\OOO)} \le M     \eqno{(1.11)}
$$
with arbitrarily chosen constant $M>0$.
\\
{\bf Case (1.2).}
\\
Then for $0<t_0 < T$, there exist constants $C>0$ and $\theta \in (0,1)$
dependent on $t_0$ and $M$ such that 
$$
\Vert u(\cdot,t_0)\Vert_{L^2(\OOO)} 
\le C(\Vert u(\cdot,T)\Vert^{\theta}_{H^1(\OOO)}
+ \Vert u(\cdot,T)\Vert_{H^1(\OOO)}).    \eqno{(1.12)}
$$
\\
{\bf Case (1.3).}
\\
We further assume that $\sigma(x)$ in (1.5) satisfies
$$
\sigma(x) > 0 \quad \mbox{for $x \in \ppp\OOO$}.      \eqno{(1.13)}
$$
Then we have (1.12).
}
\\
{\bf Theorem 2 (case $t_0=0$).}
\\
{\it
We assume (1.8) and (1.9), and $u \in H^{2,1}(Q)$ satisfy 
$\ppp_tu, \ppp_t^2u\in H^{2,1}(Q)$ and (1.1), and
$$
\sum_{k=0}^2 \Vert \ppp_t^ku(\cdot,0)\Vert_{H^1(\OOO)} \le M     \eqno{(1.14)}
$$
with arbitrarily chosen constant $M>0$.
\\
{\bf Case (1.2).}
\\
Then, for any $\alpha \in (0,1)$, there exists a constant $C>0$ such that 
$$
\Vert u(\cdot,0)\Vert_{L^2(\OOO)} \le 
C\left( \log \frac{1}{D}\right)^{-\alpha},             \eqno{(1.15)}
$$
provided that 
$$
D := \sum_{k=0}^2 \Vert \ppp_t^ku(\cdot,T)\Vert_{H^1(\OOO)}
$$
is small.
\\
{\bf Case (1.3).}
\\
Assume (1.13) additionally.  Then estimate (1.15) holds.
}
\\

Our method is quite feasible, and is applicable for example,
to semilinear equations.  We can consider more comprehensive 
class of nonlinear equations, but for discussing the essence, we are 
restricted to the following case with $0<t_0<T$:
$$
\left\{ \begin{array}{rl}
& \ppp_tu(x,t) = \sumij \ppp_i(a_{ij}(x,t)\ppp_ju(x,t))
+ c(x,t)u + f(x,t,u(x,t)), \quad (x,t) \in Q, \\
& \ppp_{\nu_A}u + r(x)u = 0 \quad \mbox{or} \quad 
u=0 \quad \mbox{on $\ppp\OOO \times (0,T)$}, 
\end{array}\right.
                                          \eqno{(1.16)}
$$
where $f(x,t,\eta)$, $x\in \ooo{\OOO}$, $0\le t \le T$ and
$\eta \in \R$, satisfies 
$$
f, \, \ppp_{\eta}f \in C(\ooo{\OOO} \times [0,T] \times \R).
                                               \eqno{(1.17)}
$$
Then we can prove
\\
{\bf Theorem 3.}
\\
{\it 
Let $0<t_0<T$.  Let $u,v \in H^{2,1}(Q)$ satisfy (1.16) and 
$$
\Vert u(\cdot,0)\Vert_{H^1(\OOO)}, \, \Vert v(\cdot,0)\Vert_{H^1(\OOO)}
\le M, \quad \Vert u\Vert_{L^{\infty}(Q)}, \, 
\Vert v\Vert_{L^{\infty}(Q)} \le M
$$
with arbitrarily chosen constant $M>0$.  Then there exist constants 
$C>0$ and $\theta \in (0,1)$, depending on $t_0$ and $M$, such that
$$
\Vert u(\cdot,t) - v(\cdot,t)\Vert_{L^2(\OOO)}
\le C(\Vert u(\cdot,T) - v(\cdot,T)\Vert^{\theta}_{H^1(\OOO)}
+ \Vert u(\cdot,T) - v(\cdot,T)\Vert_{H^1(\OOO)}).
$$
}

This article is composed of four sections.  In Section 2, we prove
a key estimate of a Carleman type.  Sections 3 and 4 are devoted to 
the proofs of Theorem 1 and 3, and Theorem 2 respectively.
%In Section 5, we provide concluding remarks.
%
%
%
%
\section{Key estimate of a Carleman estimate}

We set 
$$
\va(t) = e^{\la t}, \quad t>0,
$$
where a constant $\la>0$ is chosen later.
We state the key inequality without assumptions (1.8) and (1.9)
on the degeneracy.
\\
{\bf Lemma 1.}
\\
{\it
{\bf Case (1.2).}
\\
There exists a constant $\la_0>0$ such that for any $\la > \la_0$,
we can choose a constant $s_0(\la) > 0$ satisfying: there exist constants
$C=C(s_0,\la_0) > 0$ and $C_0 > 0$ such that  
$$
\int_Q \Biggl\{ \frac{1}{s\va} \vert \ppp_tu\vert^2
+ \sumij (\la a_{ij} - C_0\ppp_ta_{ij})(\ppp_iu)(\ppp_ju) 
+ s\la^2\va \vert u\vert^2\Biggr\} \weight dxdt
                                               \eqno{(2.1)}
$$
\begin{align*}
\le &C\int_Q \vert F\vert^2\weight dxdt
+ C\int_Q \sum_{j=1}^d \vert b_j\ppp_ju\vert^2 e^{2s\va} dxdt \\
+ & C(s\la \va(T)\Vert u(\cdot,T)\Vert^2_{L^2(\OOO)}
+ \Vert u(\cdot,T)\Vert^2_{H^1(\OOO)})e^{2s\va(T)}\\
+& C(s\la\Vert u(\cdot,0)\Vert^2_{L^2(\OOO)}
+ \Vert u(\cdot,0)\Vert^2_{H^1(\OOO)})e^{2s}
\end{align*}
for all $s > s_0$ and all $u \in H^{2,1}(Q)$ satisfying
$Lu=F$ in $Q$ and the boundary condition (1.2).
\\
{\bf Case (1.3).}
\\
There exists a constant $\la_0>0$ such that for any $\la > \la_0$,
we can choose a constant $s_0(\la) > 0$ satisfying: there exist constants
$C=C(s_0,\la_0) > 0$ and $C_0 > 0$ such that  
$$
\int_Q \Biggl\{ \frac{1}{s\va} \vert \ppp_tu\vert^2
+ \sumij (\la a_{ij} - C_0\ppp_ta_{ij})(\ppp_iu)(\ppp_ju) 
+ s\la^2\va \vert u\vert^2\Biggr\} \weight dxdt
                                    \eqno{(2.2)}
$$
\begin{align*}
\le & C\int^T_0 \int_{\ppp\OOO} \la \vert u\vert^2 e^{2s\va} dS dt\\
+ &C\int_Q \vert F\vert^2\weight dxdt
+ C\int_Q \sum_{j=1}^d \vert b_j\ppp_ju\vert^2 e^{2s\va} dxdt \\
+ & C(s\la \va(T)\Vert u(\cdot,T)\Vert^2_{L^2(\OOO)}
+ \Vert u(\cdot,T)\Vert^2_{H^1(\OOO)})e^{2s\va(T)}
+ C(s\la\Vert u(\cdot,0)\Vert^2_{L^2(\OOO)}
+ \Vert u(\cdot,0)\Vert^2_{H^1(\OOO)})e^{2s}
\end{align*}
for all $s > s_0$ and all $u \in H^{2,1}(Q)$ satisfying
$Lu=F$ in $Q$ and the boundary condition (1.3).

Here the constant $C_0$ is independent of $s_0, \la_0$.}

We emphasize that we do not assume (1.8) and (1.9), so that Lemma 1 
holds true if 
$$
\sumij a_{ij}(x,t)\xi_i\xi_j \ge 0, \quad (x,t)\in \ooo{Q},\, 
\xi_1, ..., \xi_d \in \R.
$$
In particular, 
\\
{\bf Proposition 1.}
\\
{\it 
We assume
$$
\left\{ \begin{array}{rl}
& \mbox{$a_{ij}$, $1\le i,j \le d$ are time-independent and 
$b_1 = \cdots = b_d = 0$ on $\ooo{Q}$}, \\
& \sumij a_{ij}(x)\xi_i\xi_j \ge 0, \quad x \in \ooo{\OOO},\, 
\xi_1, ..., \xi_d \in \R. \\
\end{array}\right.
$$
Then exists a constant $\la_0>0$ such that for any $\la > \la_0$,
we can choose a constant $s_0(\la) > 0$ satisfying: there exists a constant
$C=C(s_0,\la_0) > 0$ such that  
\begin{align*}
& \int_Q \Biggl\{ \frac{1}{s\va} \vert \ppp_tu\vert^2
+ s\la^2\va \vert u\vert^2\Biggr\} \weight dxdt
\le C\int_Q \vert F\vert^2\weight dxdt\\
+ &C(s\la \va(T)\Vert u(\cdot,T)\Vert^2_{L^2(\OOO)}
+ \Vert u(\cdot,T)\Vert^2_{H^1(\OOO)})e^{2s\va(T)}
+ C(s\la\Vert u(\cdot,0)\Vert^2_{L^2(\OOO)}
+ \Vert u(\cdot,0)\Vert^2_{H^1(\OOO)})e^{2s}
\end{align*}
for all $s > s_0$ and all $u \in H^{2,1}(Q)$ satisfying
$Lu=F$ in $Q$ and the boundary condition (1.2),
and
\begin{align*}
& \int_Q \Biggl\{ \frac{1}{s\va} \vert \ppp_tu\vert^2
+ s\la^2\va \vert u\vert^2\Biggr\} \weight dxdt
\le C\int^T_0 \int_{\ppp\OOO} \la \vert u\vert^2 e^{2s\va} dS dt\\
+ &C\int_Q \vert F\vert^2\weight dxdt\\
+ & C(s\la \va(T)\Vert u(\cdot,T)\Vert^2_{L^2(\OOO)}
+ \Vert u(\cdot,T)\Vert^2_{H^1(\OOO)})e^{2s\va(T)}
+ C(s\la\Vert u(\cdot,0)\Vert^2_{L^2(\OOO)}
+ \Vert u(\cdot,0)\Vert^2_{H^1(\OOO)})e^{2s}
\end{align*}
for all $s > s_0$ and all $u \in H^{2,1}(Q)$ satisfying
$Lu=F$ in $Q$ and the boundary condition (1.3).
}
\\
{\bf Remark 2.1.}
\\
As a method with a similar spirit, we can refer to the weight
energy method.  Concerning the weight energy method, there are many 
papers and see monographs Ames and Straughan \cite{AS}, 
Lees and Protter \cite{LePr}, Payne \cite{Pay},  
and the references
therein.  Except for Murray and Protter \cite{MuPr} for equations of
hyperbolic types, all the papers use just $t$ as weight function, and
do not use the second large parameter $\lambda$.
In Murray and Protter \cite{MuPr}, the weight function
$e^{s t^{\la}}$ is used to prove properties for the asymptotic
behaviour.

The essential diffference from the existing papers is  
the introduction of the second large parameter $\la > 0$.
Such a second large parameter 
is very flexible and gains a lot of possibility for better
estimates.  
\\
\vspace{0.2cm}
\\
{\bf Proof of Lemma 1.}
\\
{\bf First Step.}
\\
Set
\begin{align*}
&Lu:= \ppp_tu - \sumij \ppp_i(a_{ij}(x)\ppp_ju), \quad
G:= F(x,t) + \sum_{k=1}^d b_k(x,t)\ppp_ku(x,t), \\
& w = e^{s\va}u, \quad Pw = e^{s\va}L(e^{-s\va}w) = e^{s\va}G.
\end{align*}
Then 
$$
e^{s\va}\ppp_t(e^{-s\va}w) = \ppp_tv - s\la\va w,
$$
$$
e^{s\va}\sumij \ppp_i(a_{ij}\ppp_j(we^{-s\va}))
= \sumij \ppp_i(a_{ij}\ppp_jw),
$$
and
$$
Pw = e^{s\va}L(e^{-s\va}w) = \ppp_tw - \left( s\la\va w 
+ \sumij \ppp_i(a_{ij}\ppp_jw)\right)
= e^{s\va}G.                             
$$
We have
\begin{align*}
&\Vert e^{s\va}G\Vert^2_{L^2(Q)}\\
= &\int_Q \vert \ppp_tw\vert^2 dxdt 
+ 2\int_Q (\ppp_tw)\left(
- s\la\va w - \sumij \ppp_i(a_{ij}\ppp_jw)\right) dxdt\\
+& \int_Q \left\vert s\la\va w + \sumij \ppp_i(a_{ij}\ppp_jw)
\right\vert^2 dxdt \\
\ge &\int_Q \vert \ppp_tw\vert^2 dxdt 
+ 2\int_Q \ppp_tw\left( -\sumij \ppp_i(a_{ij}\ppp_jw)\right) dxdt\\
+ &2\int_Q (\ppp_tw)(-s\la\va)w dxdt
\end{align*}
$$
=: \int_Q \vert \ppp_tw\vert^2 dxdt + J_1 + J_2.   \eqno{(2.3)}
$$
Thus
$$
\int_Q \vert G\vert^2\weight dxdt \ge J_1+J_2                \eqno{(2.4)}
$$
and
$$
\int_Q \vert \ppp_tw\vert^2 dxdt \le \int_Q \vert G\vert^2\weight dxdt 
- J_1 - J_2.                                \eqno{(2.5)}
$$
Henceforth $C_j>0$ denote generic constants which are independent of
$s, \la$.  We assume that $s>1$ and $\la > 1$.

First we consider the boundary condition (1.2).
By noting $a_{ij} = a_{ji}$, the boundary condition (1.2) and integration by 
parts yields
$$
J_ 1 
= -2\int_Q (\ppp_tw)\sumij \ppp_i(a_{ij}\ppp_jw) dxdt 
= 2\int_Q \sumij (\ppp_i\ppp_tw)a_{ij}(\ppp_jw) dxdt 
                                                          \eqno{(2.6)}
$$
\begin{align*}
=& 2\int_Q \sum_{i>j} a_{ij}((\ppp_jw)\ppp_i\ppp_tw
+ (\ppp_iw)\ppp_j\ppp_tw) \,dxdt
+ 2\sum_{i=1}^n \int_Q a_{ii}(\ppp_iw)(\ppp_i\ppp_tw) \, dxdt\\
= & 2\int_Q \sum_{i>j} a_{ij} \ppp_t((\ppp_iw)(\ppp_jw)) \,dxdt
+ \int_Q \sum_{i=1}^n a_{ii}\ppp_t((\ppp_iw)^2) \,dxdt\\
=& \int_Q \sumij a_{ij}\ppp_t((\ppp_iw)\ppp_jw) dxdt \\
= & - \int_Q \sumij (\ppp_ta_{ij})(\ppp_iw)\ppp_jw \, dxdt
+ \sumij \left[a_{ij}(\ppp_iw)(\ppp_jw)\right]^{t=T}_{t=0} dx \\
= & -\int_Q \sumij (\ppp_ta_{ij})(\ppp_iw)(\ppp_jw)\, dxdt \\
+ & \int_{\Omega} \sumij (a_{ij}(x,T)(\ppp_iw)(x,T)(\ppp_jw)(x,T)
- a_{ij}(x,0)(\ppp_iw)(x,0)(\ppp_jw)(x,0)) dx.     
\end{align*}
On the other hand, 
\begin{align*}
&J_2 = -s\la\int_Q 2(\ppp_tw)w\va dxdt
= -s\la\int_Q \ppp_t(w^2)\va dxdt\\
=& s\la\int_Q (\ppp_t\va)w^2 dxdt
- s\la \int_{\Omega} \left[ \va w^2\right]^{t=T}_{t=0} dx
\end{align*}
$$
= s\la^2\int_Q \va w^2 dxdt
-s\la\int_{\Omega}(\va(T)\vert w(x,T)\vert^2 - \vert w(x,0)\vert^2)dx. 
                                                           \eqno{(2.7)}
$$
Hence
$$
\Vert e^{s\va}G\Vert^2_{L^2(Q)}
\ge s\la^2\int_Q \va w^2 dxdt 
- \int_Q \sumij (\ppp_ta_{ij})(\ppp_iw)(\ppp_jw) dxdt
$$
$$
- s\la \int_{\Omega} (\va(T)\vert w(x,T)\vert^2+\vert w(x,0)\vert^2)
dx
- C_1 \int_{\Omega} (\vert \nabla w(x,T)\vert^2 + 
\vert \nabla w(x,0)\vert^2) dx.                          \eqno{(2.8)}
$$
\\
{\bf Second Step.}
\\
On the right-hand side of (2.8), the term
$-\int_Q \sumij (\ppp_ta_{ij})(\ppp_iw)(\ppp_jw) dxdt$ appears, and so 
we have to estimate the integrals including $(\ppp_iw)(\ppp_jw)$.
For it, we consider $\int_Q (Pw)w \,dxdt$:
\begin{align*}
&\int_Q (Pw)w dxdt
= \int_Q (\ppp_tw)w dxdt - \int_Q s\la\va w^2 dxdt
- \int_Q \sumij \ppp_i(a_{ij}\ppp_jw)w\, dxdt\\
=: & I_1 + I_2 + I_3.
\end{align*}
We have
\begin{align*}
&\vert I_1\vert  = \left\vert 
\int_Q (\ppp_tw)w dxdt \right\vert 
= \left\vert \frac{1}{2}\int_Q \ppp_t(w^2) dxdt\right\vert\\
=& \left\vert \frac{1}{2}\int_{\Omega} \left[ \vert w(x,t)\vert^2
\right]^{t=T}_{t=0} dx \right\vert 
\le \frac{1}{2}\int_{\Omega} (\vert w(x,T)\vert^2 
+ \vert w(x,0)\vert^2) dx.
\end{align*}
Next
$$
\vert I_2\vert  = \left\vert - \int_Q s\la\va w^2 dxdt\right\vert
\le C_2\int_Q s\la\va w^2 dxdt
$$
and
$$
I_3 = -\sumij\int_Q \ppp_i(a_{ij}\ppp_jw)w dxdt
= \sumij\int_Q a_{ij}(\ppp_jw)(\ppp_iw) dxdt.
%\ge& \sigma_1 \int_Q \vert \nabla w\vert^2 dxdt.        
                                       \eqno{(2.9)}
$$
Hence
$$
\int_Q \la(Pw)w\,dxdt
\ge \la\int_Q \sumij a_{ij}(\ppp_iw)(\ppp_jw) dxdt 
- C_2\int_Q s\la^2\va w^2 dxdt
$$
$$
- \frac{1}{2}\la\int_{\Omega} (\vert w(x,T)\vert^2 
+ \vert w(x,0)\vert^2) dx.                          \eqno{(2.10)}
$$
On the other hand,
$$
\left\vert \int_Q \la (Pw)w \, dxdt \right\vert 
\le \Vert Pw\Vert_{L^2(Q)}(\la\Vert w\Vert_{L^2(Q)})
\le \frac{1}{2}\Vert Pw\Vert^2_{L^2(Q)}
+ \frac{\la^2}{2}\Vert w\Vert_{L^2(Q)}^2
$$
$$
= \frac{1}{2}\int \vert G\vert^2 \weight dxdt 
+ \frac{\la^2}{2}\Vert w\Vert_{L^2(Q)}^2.     \eqno{(2.11)}
$$
Hence (2.10) yields
\begin{align*}
& \la\int_Q \sumij a_{ij}(\ppp_iw)(\ppp_jw) dxdt
\le C_2\int_Q s\la^2\va w^2 dxdt\\
+ & \frac{1}{2}\int_Q \vert Ge^{s\va}\vert^2 dxdt
+ \frac{\la^2}{2} \int_Q w^2 dxdt
+ \frac{1}{2}\la\int_{\Omega} (\vert w(x,T)\vert^2 
+ \vert w(x,0)\vert^2) dx.
\end{align*}
Estimating the first term on the right-hand side by (2.8), we obtain
\begin{align*}
&\la \int_Q \sumij a_{ij}(\ppp_iw)(\ppp_jw) dxdt\\
\le &C_3\int_Q \vert Ge^{s\va}\vert^2 dxdt
+ C_3\int_Q \sumij (\ppp_ta_{ij})(\ppp_iw)(\ppp_jw) \, dxdt 
+ C_3\int_Q \la^2w^2 dxdt\\
+& C_3\la(\Vert w(\cdot,T)\Vert^2_{L^2(\Omega)} 
+ \Vert w(\cdot,0)\Vert^2_{L^2(\Omega)})
+ C_3(\Vert \nabla w(\cdot,T)\Vert^2_{L^2(\Omega)} 
+ \Vert \nabla w(\cdot,0)\Vert^2_{L^2(\Omega)})\\
+ &C_3s\la (\va(T)\Vert w(\cdot,T)\Vert^2_{L^2(\Omega)}
+ \Vert w(\cdot,0)\Vert^2_{L^2(\Omega)})\\
\le& C_3\int_Q \vert Ge^{s\va}\vert^2 dxdt
+ C_3\int_Q \sumij (\ppp_ta_{ij})(\ppp_iw)(\ppp_jw) \, dxdt 
+ C_3\int_Q \la^2w^2 dxdt
\end{align*}
$$
+ C_3s\la (\va(T)\Vert w(\cdot,T)\Vert^2_{L^2(\Omega)} 
+ \Vert w(\cdot,0)\Vert^2_{L^2(\Omega)})
+ C_3(\Vert \nabla w(\cdot,T)\Vert^2_{L^2(\Omega)} 
+ \Vert \nabla w(\cdot,0)\Vert^2_{L^2(\Omega)})
                                                \eqno{(2.12)}
$$
\\
{\bf Third Step.}
\\
Adding (2.8) and (2.12), we have
\begin{align*}
&\int_Q s\la^2\va w^2 dxdt 
+ \la\int_Q \sumij a_{ij}(\ppp_iw)(\ppp_jw) dxdt\\
\le& C_4\int_Q \vert Ge^{s\va}\vert^2 dxdt
+ C_3\int_Q \sumij (\ppp_ta_{ij})(\ppp_iw)(\ppp_jw) dxdt
+ C_4\int_Q \la^2w^2 dxdt\\
+& C_4(\Vert \nabla w(\cdot,T)\Vert^2_{L^2(\Omega)} 
+ \Vert \nabla w(\cdot,0)\Vert^2_{L^2(\Omega)})
+ C_4s\la (\va(T)\Vert w(\cdot,T)\Vert^2_{L^2(\Omega)}
+ \Vert w(\cdot,0)\Vert^2_{L^2(\Omega)}).                
\end{align*}
By $\va = e^{\la t} \ge 1$, we take $s>0$ and $\la > 0$ 
large to absorb the third term on the right-hand side 
into the left-hand side.
Hence, 
\begin{align*}
&\int_Q s\la^2\va w^2 dxdt 
+ \int_Q \sumij (\la a_{ij} - C_3\ppp_ta_{ij})(\ppp_iw)(\ppp_jw) dxdt \\
\le &C_5\int_Q \vert Ge^{s\va}\vert^2 dxdt\\
+ & C_5s\la(\va(T)\Vert w(\cdot,T)\Vert^2_{L^2(\Omega)}
+ \Vert w(\cdot,0)\Vert^2_{L^2(\Omega)})
\end{align*}
$$
+ C_5(\Vert \nabla w(\cdot,T)\Vert^2_{L^2(\Omega)}
+ \Vert \nabla w(\cdot,0)\Vert^2_{L^2(\Omega)}).    
                                                    \eqno{(2.13)}
$$
\\
{\bf Fourth Step.}
\\
Next we will estimate $\vert \ppp_tw\vert^2$.  Since 
$u = e^{-s\va}w$, we have $\ppp_tu = -s\la\va e^{-s\va}w
+ e^{-s\va}\ppp_tw$, and
$$
\frac{1}{s\va}\vert \ppp_tu\vert^2\weight
\le 2s\la^2\va w^2 + \frac{2}{s\va}\vert \ppp_tw\vert^2.
$$
Let $\ep \in \left(0, \frac{1}{2}\right)$ be a constant 
which we choose later.  We note that  
$\frac{1}{s\va} = \frac{1}{se^{\la t}} \le \frac{1}{s} \le 1$
for $s\ge 1$.
Therefore, for all large $s > 0$ and $\la>0$, we have
\begin{align*}
&\int_Q \frac{\ep}{s\va}\vert \ppp_tu\vert^2\weight dxdt
\le \int_Q 2\ep s \la^2\va w^2 dxdt 
+ \int_Q \frac{2\ep}{s\va}\vert \ppp_tw\vert^2 dxdt\\
\le& 2\ep\int_Q s\la^2\va w^2 dxdt 
+ \ep\int_Q \vert \ppp_tw\vert^2 dxdt
\end{align*}
$$
\le 2\ep\int_Q s\la^2\va w^2 dxdt + \ep\int_Q G^2\weight dxdt
+ \ep (-J_1 - J_2)                  \eqno{(2.14)}
$$
by (2.5).  

By (2.6) and (2.7), we have
\begin{align*}
& \ep(-J_1-J_2)
= \int_Q \ep\sumij (\ppp_ta_{ij})(\ppp_iw)(\ppp_jw) dxdt 
- s\la^2\ep \int_Q \va w^2 dxdt\\
-& \ep\int_{\OOO} \sumij (a_{ij}(x,T)(\ppp_iw)(x,T)(\ppp_jw)(x,T)
- a_{ij}(x,0)(\ppp_iw)(x,0)(\ppp_jw)(x,0)) dx\\
+& \ep s\la \int_{\OOO} (\va(T)\vert w(x,T)\vert^2
- \vert w(x,0)\vert^2) dx,
\end{align*}
and the substitution of this into (2.13) yields
\begin{align*}
& \ep\int_Q \frac{1}{s\va} \vert \ppp_tu\vert^2 \weight dxdt
\le 2\ep\int_Q s\la^2\va w^2 dxdt + \ep\int_Q \vert G\vert^2\weight dxdt\\
+& \sumij \int_Q \ep(\ppp_ta_{ij})(\ppp_iw)(\ppp_jw) dxdt
- \ep\int_Q s\la^2\va w^2 dxdt\\
+& C_6\ep\int_{\OOO} (\vert \nabla w(x,T)\vert^2
+ \vert \nabla w(x,0)\vert^2) dx
\end{align*}
$$
+ \ep s \la \int_{\OOO} (\va(T)\vert w(x,T)\vert^2
+ \vert w(x,0)\vert^2) dx.             \eqno{(2.15)}
$$
Adding (2.13) and (2.15), we obtain
\begin{align*}
& \int_Q s\la^2\va w^2 dxdt
+ \int_Q \sumij (\la a_{ij} - C_3\ppp_ta_{ij})
(\ppp_iw)(\ppp_jw) dxdt 
+ \ep\int_Q \frac{1}{s\va} \vert \ppp_tu\vert^2 \weight dxdt\\
\le& C_7\int_Q \vert G\vert^2 \weight dxdt 
+ \ep\int_Q s\la^2 \va w^2 dxdt \\
+& \int_Q \ep \sumij (\ppp_ta_{ij})(\ppp_iw)(\ppp_jw) dxdt\\
+& C_7(s\la \va(T)\Vert w(\cdot,T)\Vert^2_{L^2(\OOO)}
+ \Vert \nabla w(\cdot,T)\Vert_{L^2(\OOO)}^2)\\
+& C_7(s\la \Vert w(\cdot,0)\Vert^2_{L^2(\OOO)}
+ \Vert \nabla w(\cdot,0)\Vert_{L^2(\OOO)}^2).
\end{align*}
Choosing $0 < \ep < \frac{1}{2}$, we can absorb the second term on the 
right-hand side into the left-hand side, we obtain
\begin{align*}
& \ep\int_Q \frac{1}{s\va} \vert \ppp_tu\vert^2 \weight dxdt\\
+ & \int_Q \sumij (\la a_{ij} - (C_3+\ep)\ppp_ta_{ij})
(\ppp_iu)(\ppp_ju) \weight dxdt 
+ \frac{1}{2}\int_Q s\la^2\va w^2 dxdt\\
\le& C_7\int_Q \vert G\vert^2 \weight dxdt \\
+ & C_7(s\la \va(T)\Vert w(\cdot,T)\Vert^2_{L^2(\OOO)}
+ \Vert \nabla w(\cdot,T)\Vert_{L^2(\OOO)}^2)\\
+& C_7(s\la \Vert w(\cdot,0)\Vert^2_{L^2(\OOO)}
+ \Vert \nabla w(\cdot,0)\Vert_{L^2(\OOO)}^2).
\end{align*}
Thus the proof of Lemma 1 is complete in the case (1.2).
\\
{\bf Fifth Step.}
\\
We will prove in the case (1.3).
We recall the additional assumption (1.13).  
In the above arguments, we need to modify only (2.6) and (2.9).
That is,
\begin{align*}
&J_1 = -2\int_Q (\ppp_tw)\sumij \ppp_i(a_{ij}\ppp_jw) dxdt\\
= &2\int_Q \sumij (\ppp_i\ppp_tw)a_{ij}(\ppp_jw) dxdt
- 2\int_{\ppp\Omega\times (0,T)}
\sumij a_{ij}(\ppp_jw)\nu_i(\ppp_tw) dS dt\\
=& 2\int_Q \sumij a_{ij}(\ppp_i\ppp_tw)(\ppp_jw) dxdt
- 2\int_{\ppp\Omega\times(0,T)} (\ppp_{\nu_A}w)(\ppp_tw) dS dt.
\end{align*}
By (1.3), we obtain
\begin{align*}
&- 2\int_{\ppp\Omega\times(0,T)} (\ppp_{\nu_A}w)(\ppp_tw)
dSdt
= 2\int_{\ppp\Omega\times(0,T)}rw(\ppp_tw)dSdt\\
=& \int_{\ppp\Omega\times(0,T)} r\ppp_t(w^2) dSdt
= \int_{\ppp\Omega} \left[ rw^2\right]^{t=T}_{t=0} dS.
\end{align*}
Thus, similarly to (2.6), we have
\begin{align*}
& -J_1 = -2\int_Q \sumij a_{ij}(\ppp_i\ppp_tw)(\ppp_jw) dxdt
-\int_{\ppp\Omega} \left[ rw^2\right]^{t=T}_{t=0} dS\\
=& \int_Q \sumij (\ppp_ta_{ij})(\ppp_iw)(\ppp_jw) dxdt \\
- & \int_Q \sumij (a_{ij}(x,T)(\ppp_iw)(x,T)(\ppp_jw)(x,T) 
- a_{ij}(x,0)(\ppp_iw)(x,0)\ppp_jw(x,0)) dx\\
-& \int_{\OOO} r(w^2(x,T) - w^2(x,0)) dS
\end{align*} 
and by the trace theorem: $\Vert u\Vert_{L^2(\ppp\Omega)}
\le C\Vert u\Vert_{H^1(\Omega)}$, we can obtain the same estimate
$$
-J_1 \le \int_Q \sumij (\ppp_ta_{ij})(\ppp_iw)(\ppp_jw) dxdt 
+ C_8(\Vert w(\cdot,T)\Vert_{H^1(\Omega)}^2
+ \Vert w(\cdot,0)\Vert_{H^1(\Omega)}^2).         \eqno{(2.16)}
$$
Therefore (2.4) and (2.7) yield
\begin{align*}
& \int_Q \vert G\vert^2 \weight dxdt \ge J_1 + J_2\\
\ge& \int_Q s\la^2\va w^2 dxdt 
- s\la(\va(T)\Vert w(\cdot,T)\Vert^2_{L^2(\OOO)}
+ \Vert w(\cdot,0)\Vert^2_{L^2(\OOO)})\\
- & \int_Q \sumij (\ppp_ta_{ij})(\ppp_iw)(\ppp_jw) dxdt 
- C_8(\Vert w(\cdot,T)\Vert_{H^1(\Omega)}^2
+ \Vert w(\cdot,0)\Vert_{H^1(\Omega)}^2),
\end{align*}
that is,
$$
\int_Q s\la^2\va w^2 dxdt              \eqno{(2.17)}
$$
\begin{align*}
\le& \int_Q \vert G\vert^2 \weight dxdt
+ s\la(\va(T)\Vert w(\cdot,T)\Vert^2_{L^2(\OOO)}
+ \Vert w(\cdot,0)\Vert^2_{L^2(\OOO)})\\
+& \int_Q \sumij (\ppp_ta_{ij})(\ppp_iw)(\ppp_jw) dxdt 
+ C_8(\Vert w(\cdot,T)\Vert_{H^1(\Omega)}^2
+ \Vert w(\cdot,0)\Vert_{H^1(\Omega)}^2).
\end{align*}

Moreover, for $I_3$ in (2.9), by (1.5) we have
\begin{align*}
&I_3 = \sumij \int_Q a_{ij}(\ppp_iw)(\ppp_jw) dx
- \int_{\ppp\Omega\times(0,T)} w(\ppp_{\nu_A}w)\, dSdt\\
=& \sumij \int_Q a_{ij}(\ppp_iw)(\ppp_jw) dx
+ \int_{\ppp\Omega\times(0,T)} rw^2 dSdt\\
\ge& \int_Q \sumij a_{ij}(\ppp_iw)(\ppp_jw) dxdt
- C_9\Vert w\Vert^2_{L^2(0,T;L^2(\ppp\Omega))}.
\end{align*}
Hence, since we can obtain the same estimates for $I_1$ and $I_2$, we see
$$
\int_Q \la(Pw)w dxdt = \la I_1 + \la I_2 + \la I_3
                                     \eqno{(2.18)}
$$
\begin{align*}
\ge& \int_Q \sumij \la a_{ij}(\ppp_iw)(\ppp_jw) dxdt 
- C_{10}\Vert w\Vert^2_{L^2(0,T;L^2(\ppp\OOO))}\\
- & \frac{1}{2}\la(\Vert w(\cdot,T)\Vert^2_{L^2(\OOO)}
+ \Vert w(\cdot,0)\Vert^2_{L^2(\OOO)})
- C_{10}\int_Q s\la^2\va w^2 dxdt.
\end{align*}

Similarly to (2.12), using (2.11), by (2.18) we obtain
\begin{align*}
&\int_Q \sumij \la a_{ij}(\ppp_iw)(\ppp_jw) dxdt\\
\le& \int_Q \la(Pw)w dxdt 
+ C_{10}\la \Vert w\Vert^2_{L^2(0,T;L^2(\ppp\OOO))}
+ \frac{1}{2}\la(\Vert w(\cdot,T)\Vert^2_{L^2(\OOO)}
+ \Vert w(\cdot,0)\Vert^2_{L^2(\OOO)})\\
+ & C_{10}\int_Q s\la^2\va w^2 dxdt\\
\le& \frac{1}{2} \Vert Ge^{s\va}\Vert^2_{L^2(Q)}
+ \frac{\la^2}{2}\Vert w\Vert^2_{L^2(Q)}
+ C_{10}\la \Vert w\Vert^2_{L^2(0,T;L^2(\ppp\OOO))} \\
+ & \frac{1}{2}\la(\Vert w(\cdot,T)\Vert^2_{L^2(\OOO)}
+ \Vert w(\cdot,0)\Vert^2_{L^2(\OOO)}
+ C_{10}\int_Q s\la^2\va w^2 dxdt.
\end{align*}
Hence (2.17) implies 
$$
\int_Q \sumij \la a_{ij}(\ppp_iw)(\ppp_jw) dxdt
                                                \eqno{(2.19)}
$$
\begin{align*}
\le& \frac{1}{2} \Vert Ge^{s\va}\Vert^2_{L^2(Q)}
+ C_{11}\la^2 \Vert w\Vert^2_{L^2(Q)}\\
+& C_{11}\la \Vert w\Vert^2_{L^2(0,T;L^2(\ppp\OOO))}
+ \frac{1}{2}\la( \Vert w(\cdot,T)\Vert^2_{L^2(\OOO)}
+ \Vert w(\cdot,0)\Vert^2_{L^2(\OOO)})\\
+& C_{11}\int_Q \vert G\vert^2 \weight dxdr
+ C_{11}s\la(\va(T)\Vert w(\cdot,T)\Vert^2_{L^2(\OOO)}
+ \Vert w(\cdot,0)\Vert^2_{L^2(\OOO)})\\
+& \int_Q \sumij (C_{11}\ppp_ta_{ij})(\ppp_iw)(\ppp_jw) dxdt\\
+ &C_{11}(\Vert w(\cdot,T)\Vert_{H^1(\Omega)}^2
+ \Vert w(\cdot,0)\Vert_{H^1(\Omega)}^2).
\end{align*}
Thus (2.17) and (2.19) yield
\begin{align*}
& \int_Q \sumij (\la a_{ij} - C_{11}\ppp_ta_{ij})(\ppp_iw)(\ppp_jw) dxdt
+ \int_Q s\la^2\va w^2 dxdt\\
\le& C_{12}\int_Q \vert G\vert^2 \weight dxdt 
+ C_{12}\la \Vert w\Vert^2_{L^2(0,T;L^2(\ppp\OOO))}
+ C_{12}(\Vert w(\cdot,T)\Vert^2_{H^1(\OOO)}
+ \Vert w(\cdot,0)\Vert^2_{H^1(\OOO)})
\end{align*}
$$
+ C_{12}s\la (\va(T)\Vert w(\cdot,T)\Vert^2_{L^2(\OOO)}
+ \Vert w(\cdot, 0)\Vert^2_{L^2(\OOO)}).              \eqno{(2.20)}
$$
Here we used 
$$
\frac{1}{2}(\Vert w(\cdot,T)\Vert^2_{L^2(\OOO)}
+ \Vert w(\cdot,0)\Vert^2_{L^2(\OOO)})
\le C_{12}s\la (\va(T)\Vert w(\cdot,T)\Vert^2_{L^2(\OOO)}
+ \Vert w(\cdot,0)\Vert^2_{L^2(\OOO)})
$$
and we absorbed the term
$$
C_{11}\la^2\Vert w\Vert^2_{L^2(Q)} = C_{11}\int_Q \la^2w^2 dxdt
$$
into $\int_Q s\la^2\va w^2 dxdt$.

Next as for the estimate of $\ppp_tw$, we can proceed
similarly to the argument starting (2.14) as follows.
In terms of (2.7) and (2.16), we have 
$$
 \ep(-J_1 - J_2)                    \eqno{(2.21)}
$$
\begin{align*}
\le& C_{13}\ep\int_Q \sumij (\ppp_t a_{ij})(\ppp_iw)(\ppp_jw) dxdt
+ C_{13}\ep(\Vert w(\cdot,T)\Vert^2_{H^1(\OOO)}
+ \Vert w(\cdot,0)\Vert^2_{H^1(\OOO)})\\
- &C_{13}\ep\int_Q s\la^2 \va w^2 dxdt
+ \ep s\la (\va(T)\Vert w(\cdot,T)\Vert^2_{L^2(\OOO)}
+ \Vert w(\cdot,0)\Vert^2_{L^2(\OOO)}).
\end{align*}
Thus, by (2.21) and (2.14), we can estimate 
$\int_Q \frac{1}{s\va}\vert \ppp_tu\vert^2
\weight dxdt$, and so (2.21) completes the proof of 
Lemma 1 in the case (1.3).
$\blacksquare$
\section{Proofs of Theorems 1 and 3}

In Lemma 1, in order to estimate $\left( \frac{1}{s\va}\vert \ppp_tu\vert^2
+ s\la^2\va \vert u \vert^2\right) e^{2s\va}$, we have to assume 
(1.8) and (1.9) because the term
$$
\sumij (\la a_{ij} - C_0\ppp_ta_{ij})(\ppp_iu)(\ppp_ju)\weight
$$
is not necessarily non-negative.
More precisely, we can state 
\\
{\bf Lemma 2.}
\\
{\it
We assume (1.8) and (1.9).   In case (1.3), we further assume
(1.13).
Then there exists a constant $\la_0>0$ such that for any $\la > \la_0$,
we can choose a constant $s_0(\la) > 0$ satisfying: there exists a constant
$C=C(s_0,\la_0) > 0$ such that  
$$
\int_Q \Biggl\{ \frac{1}{s\va}\vert \ppp_tu\vert^2
+ s\la^2\va \vert u\vert^2\Biggr\} \weight dxdt
                                    \eqno{(3.1)}
$$
\begin{align*}
\le &C\int_Q \vert F\vert^2\weight dxdt
+ C(s\la \va(T)\Vert u(\cdot,T)\Vert^2_{L^2(\OOO)}
+ \Vert u(\cdot,T)\Vert^2_{H^1(\OOO)})e^{2s\va(T)}\\
+& C(s\la\Vert u(\cdot,0)\Vert^2_{L^2(\OOO)}
+ \Vert u(\cdot,0)\Vert^2_{H^1(\OOO)})e^{2s}
\end{align*}
for all $s > s_0$ whenever a solution $u \in H^{2,1}(Q)$ satisfies 
(1.1) with (1.2) or (1.3).
}
\\
{\bf Proof of Lemma 2.}
\\
{\bf Case (1.2).}
\\
First let $b:= (b_1,  \cdots, b_d) = 0$ in $\OOO$.
By (1.8), we can choose $\la>0$ large such that 
$$
\sumij (\la a_{ij} - C_0\ppp_ta_{ij})\xi_i\xi_j \ge 0
\quad \mbox{for $(x,t) \in \ooo{Q}$ and 
$\xi_1, ..., \xi_d \in \R$.}                     \eqno{(3.2)}
$$
Indeed, by (1.8) and (1.5) with $\sigma(x) \ge 0$, we verify that 
\begin{align*}
& \sumij (\la a_{ij} - C_0\ppp_ta_{ij})(\ppp_iw)(\ppp_jw)\\
=& \sumij (\la_1C_0 a_{ij} - C_0\ppp_ta_{ij})(\ppp_iw)(\ppp_jw)
+ (\la-\la_1C_0)\sumij a_{ij}(\ppp_iw)(\ppp_jw)\\
\ge& C_0\sumij (\la_1 a_{ij} - \ppp_ta_{ij})(\ppp_iw)(\ppp_jw)\ge 0
\quad \mbox{on $\ooo{Q}$},
\end{align*}
provided that $\la>0$ is so large that $\la - \la_1C_0 \ge 0$, which 
verifies (3.2).
\\

We set $\va(t) = e^{\la t}$.  Since $b=0$ in $\OOO$, 
by Lemma 1, we can readily obtain (3.1).

Second let $b \not\equiv 0$ in $\OOO$.  
Then (1.8) and (1.5) yield
\begin{align*}
& \sumij (\la a_{ij} - C_0\ppp_ta_{ij})(\ppp_iw)(\ppp_jw)\\
=& \sumij (\la- \la_1C_0) a_{ij}(\ppp_iw)(\ppp_jw)
+ C_0\sumij (\la_1a_{ij}- \ppp_ta_{ij})(\ppp_iw)(\ppp_jw)\\
\ge& \sumij (\la - \la_1C_0)a_{ij}(\ppp_iw)(\ppp_jw)
\ge (\la - \la_1C_0)\sigma(x)\vert \nabla w(x,t)\vert^2,
\quad (x,t) \in \ooo{Q}.
\end{align*}
Therefore, (1.9) implies
$$
\int_Q \sumij (\la a_{ij} - C_0\ppp_ta_{ij})(\ppp_iw)(\ppp_jw) dxdt
- C\int_Q \sum_{j=1}^d \vert b_j\ppp_jw\vert^2 dxdt
                                                \eqno{(3.3)}
$$
\begin{align*}
\ge & \int_Q (\la - \la_1C_0)\sigma(x)\vert \nabla w(x,t)\vert^2 dxdt
- C\int_Q \sigma(x)\vert \nabla w(x,t)\vert^2 dxdt\\
\ge &(\la - \la_1C_0 - C)\int_Q \sigma(x)\vert \nabla w(x,t)\vert^2 dxdt.
\end{align*}
Hence, choosing $\la > 0$ sufficiently large, we can absorb the second term
on the right-hand side of (2.1) into the left-hand side.

Fixing such sufficiently large $\la>0$, we obtain (3.1). 
\\
{\bf Case (1.3).}
\\
We assume that as a neighborhood of $\ppp\OOO$, by means of (1.13), we can 
find a subdomain $\OOO'$ satisfying
$\ooo{\OOO'} \subset \OOO$ and a constant $\ep_1 > 0$ such that 
$$
\ppp\OOO' \supset \ppp\OOO, \quad \sigma(x) \ge \ep_1 \quad 
\mbox{for $x \in \OOO'$}.                                     \eqno{(3.4)}
$$

We fix $\delta > 0$ sufficiently small.  Then, 
the interpolation inequality and the trace theorem (e.g., Adams \cite{Ad})
imply that for any $\delta > 0$ there exists a constant $C_2(\delta) > 0$
such that 
$$
\Vert u(\cdot,t)\Vert_{L^2(\ppp\OOO)} 
\le C_1\Vert u(\cdot,t)\Vert_{H^{\hhalf+\delta_0}(\OOO')}
\le \delta\Vert\nabla u(\cdot,t)\Vert_{L^2(\OOO')}
+ C_2(\delta)\Vert u(\cdot,t)\Vert_{L^2(\OOO')}.      \eqno{(3.5)}
$$
Since (3.3) holds true also in the case (1.3), from (2.2) in Lemma 1,
we can derive
$$
\int_Q \left( \frac{1}{s\va}\vert \ppp_tu\vert^2
+ s\la^2\va \vert u\vert^2\right)) \weight dxdt
+ (\la-\la_1C_0-C)\int_Q \sigma(x)\vert \nabla u\vert^2 \weight
dxdt                                                   \eqno{(3.6)}
$$
\begin{align*}      
\le& C\int^T_0 \int_{\ppp\OOO} \la \vert u\vert^2 \weight dS dt
+ C\int_Q \vert F\vert^2 \weight dxdt \\
+& C(s\la \va(T)\Vert u(\cdot,T)\Vert^2_{L^2(\OOO)}
+ \Vert u(\cdot,T)\Vert^2_{H^1(\OOO)})e^{2s\va(T)}\\
+& C(s\la\Vert u(\cdot,0)\Vert^2_{L^2(\OOO)}
+ \Vert u(\cdot,0)\Vert^2_{H^1(\OOO)})e^{2s}.
\end{align*}
Choosing $\la > 0$ further large such that $\la - \la_1C_0 - C > 0$,
we see by (3.3) that 
$$
(\la-\la_1C_0-C)\int_Q \sigma(x)\vert \nabla u\vert^2 \weight dxdt
\ge (\la-\la_1C_0-C)\int_{\OOO'\times (0,T)} 
\sigma(x)\vert \nabla u\vert^2 \weight dxdt      \eqno{(3.7)}
$$
$$
\ge C_3\ep_1 \int_{\OOO'\times (0,T)} \vert \nabla u\vert^2 \weight dxdt.
$$
Therefore, in terms of (3.5), we obtain
\begin{align*}
& \la \Vert u(\cdot,t)e^{s\va(t)}\Vert^2_{L^2(\ppp\OOO)}\\
\le& C\la\delta \Vert \nabla u(\cdot,t)e^{s\va(t)}\Vert^2_{L^2(\OOO')}
+ C\la C_2(\delta) \Vert u(\cdot,t)e^{s\va(t)}\Vert^2_{L^2(\OOO')}.
\end{align*}
We substitute these inequalities into the second term on the left-hand side 
and the first term on the right-hand side of (3.6) to reach
\begin{align*}
& \int_Q \left( \frac{1}{s\va}\vert \ppp_tu\vert^2
+ s\la^2\va \vert u\vert^2\right)) \weight dxdt
+ C_3\ep_1\int_{\OOO'\times (0,T)} \vert \nabla u\vert^2 \weight\\
\le& C\la\delta \int_{\OOO'\times (0,T)} \vert \nabla u\vert^2 \weight
dxdt
+ C\la C_2(\delta) \int_{\OOO'\times (0,T)} \vert u\vert^2 \weight dxdt\\
+ & C\int_Q \vert F\vert^2 \weight dxdt \\
+& C(s\la \va(T)\Vert u(\cdot,T)\Vert^2_{L^2(\OOO)}
+ \Vert u(\cdot,T)\Vert^2_{H^1(\OOO)})e^{2s\va(T)}\\
+& C(s\la\Vert u(\cdot,0)\Vert^2_{L^2(\OOO)}
+ \Vert u(\cdot,0)\Vert^2_{H^1(\OOO)})e^{2s}.
\end{align*}
Choosing $\delta > 0$ sufficiently small and $s>0$ sufficiently large, 
we can absorb the first and the second terms on the 
right-hand side into the second and the first terms on the left-hand side.
Thus the proof of Lemma 2 is complete.
$\blacksquare$
\\

Now we proceed to 
\\
{\bf Proof of Theorem 1.}
\\
Using (1.11) and $\va(t_0) \le \va(t)$ for $t_0\le t \le T$, we shrink
the integral region $Q$ on the right-hand side of (3.1), we obtain
\begin{align*}
& e^{2s\va(t_0)} \int_{\OOO\times (t_0,T)}
 \left( \frac{1}{s\va}\vert \ppp_tu\vert^2
+ s\la^2\va \vert u\vert^2\right) dxdt\\
\le& Cs\la\va(T) \Vert u(\cdot,T)\Vert^2_{H^1(\OOO)}e^{2s\va(T)}
+ Cs\la M^2e^{2s}.
\end{align*}
for $s \ge s_0$.
Fixing $\la > 0$ sufficiently large, we do not need specify the 
$\la$-dependency, and so  
$$
\int_{\OOO\times (t_0,T)}
 \left( \frac{1}{s}\vert \ppp_tu\vert^2
+ s\vert u\vert^2\right) dxdt
\le Cs D_0^2 e^{2s(\va(T)-\va(t_0))}
+ Cs M^2e^{-2s\mu(t_0)}.
$$
for $s \ge s_0$.
Here the constant $C>0$ depends on $T$ and $\la$.
We set 
$$
\mu(t_0) := \va(t_0) - 1 = e^{\la t_0} - 1 > 0, \quad
D_0:= \Vert u(\cdot,T)\Vert_{H^1(\OOO)}.
$$
Thus
$$
\Vert \ppp_tu\Vert^2_{L^2(t_0,T;L^2(\OOO))}
\le Cs^2D_0^2e^{2s(\va(T)-\va(t_0))} + Cs^2M^2e^{-2s\mu(t_0)}
                                                          \eqno{(3.8)}
$$
for all large $s>0$.
 
We note that the generic contants $C>0$, $C_j>0$ are independent of also
$t_0 \in [0,T]$, but dependent on $T$, $\OOO$, $\la$.

Since
$$
u(x,t_0) = \int^{t_0}_T \ppp_tu(x,t) dt + u(x,T), \quad x\in \OOO,
$$
we can choose a constant $C_4>0$ such that 
$$
\Vert u(\cdot,t_0)\Vert^2_{L^2(\OOO)}
\le C_4\Vert \ppp_tu\Vert_{L^2(\OOO\times (t_0,T))}^2
+ C_4\Vert u(\cdot,T)\Vert^2_{L^2(\OOO)}                    \eqno{(3.9)}
$$
for all $t_0 \in [0,T]$.
Substituting (3.8) into (3.9) and using $\va(T) > 1$, we obtain
\begin{align*}
& \Vert u(\cdot,t_0)\Vert_{L^2(\OOO)}^2
\le C_4(Cs^2D_0^2e^{2s\va(T)} + Cs^2M^2e^{-2s\mu(t_0)})
+ C_4\Vert u(\cdot,T)\Vert^2_{L^2(\OOO)}\\
\le & C_5s^2D_0^2e^{2s\va(T)} + C_5s^2M^2e^{-2s\mu(t_0)}
\le C_6D_0^2e^{3s\va(T)} + C_6M^2e^{-s\mu(t_0)}
\end{align*}
for all large $s>s_0$ and all $t_0\in [0,T]$.

Setting $C_7:= C_6e^{3s_0\va(T)}$, we can have
$$
\Vert u(\cdot,t_0)\Vert_{L^2(\OOO)}^2
\le C_7D_0^2e^{3s\va(T)} + C_7M^2e^{-s\mu(t_0)}     \eqno{(3.10)}
$$
for all large $s>0$ and all $t_0\in [0,T]$.

Now we choose $s>0$ in order to minimize the right-hand side of (3.10).
\\
{\bf Case 1: $M\le D_0$}.
\\
By setting $s=0$, the inequality (3.10) immediately yields 
$\Vert u(\cdot,t_0)\Vert_{L^2(\OOO)}^2
\le 2C_7D^2_0$.
\\
{\bf Case 2: $M>D_0$}.
\\
We choose $s>0$ such that 
$$
D_0^2e^{3s\va(T)} = M^2e^{-s\mu(t_0)}, 
$$
that is,
$$
s = \frac{2}{3\va(T)+\mu(t_0)}\log \frac{M}{D_0} > 0,
$$
so that 
$$
\Vert u(\cdot,t_0)\Vert^2_{L^2(\OOO)}
\le 2C_7M^{\frac{6\va(T)}{3\va(T)+\mu(t_0)}}
D_0^{\frac{2\mu(t_0)}{3\va(T)+\mu(t_0)}}.
$$
Setting $\theta := \frac{\mu(t_0)}{3\va(T)+\mu(t_0)}
\in (0,1)$, we can obtain (1.12).
Also in the case of (1.3), thanks to Lemma 2, the same 
argument completes the proof of Theorem 1.
$\blacksquare$
\\

Nest we provide
\\
{\bf Proof of Theorem 3.}
\\
We set $y:= u-v$.  Then
$$
\left\{ \begin{array}{rl}
& \ppp_ty = \sumij \ppp_i(a_{ij}\ppp_jy) + c(x,t)y
+ (f(x,t,u(x,t)) - f(x,t,v(x,t)), \quad (x,t) \in Q, \\
& \ppp_{\nu_A}y + ry = 0 \quad \mbox{or} \quad
y=0 \quad \mbox{on $\ppp\OOO \times (0,T)$}.
\end{array}\right.
$$
We apply Lemma 2 to obtain
$$
\int_Q \left( \frac{1}{s\va}\vert \ppp_ty\vert^2
+ s\la^2\va \vert y\vert^2\right) \weight dxdt
                                              \eqno{(3.11)}
$$
\begin{align*}
\le & \int_Q \vert f(x,t,u(x,t)) - f(x,t,v(x,t))\vert^2 \weight dxdt\\
+& Cs\la \va(T)\Vert y(\cdot,T)\Vert^2_{H^1(\OOO)} e^{2s\va(T)}
+ Cs\la\Vert y(\cdot,0)\Vert^2_{H^1(\OOO)}e^{2s}
\end{align*}
for all large $s>0$.

By (1.2) and $\Vert u\Vert_{L^{\infty}(Q)}$,
$\Vert u\Vert_{L^{\infty}(Q)} \le M$, in terms of the mean value 
theorem, we have
$$
\vert f(x,t,u(x,t)) - f(x,t,v(x,t))\vert 
= \vert f(x,t,\eta_{x,t})\vert \vert u(x,t) - v(x,t)\vert
\le C_8\vert u(x,t) - v(x,t)\vert,
$$
where $\eta_{x,t}$ is a constant depending on $x,t$ such that 
$\vert \eta_{x,t}\vert \le M$.  Therefore, we can absorb the first term
on the right-hand side of (3.11) into the left-hand side by choosing
$s>0$ sufficiently large, and we obtain
\begin{align*}
& \int_Q \left( \frac{1}{s\va}\vert \ppp_ty\vert^2
+ s\la^2\va \vert y\vert^2\right) \weight dxdt  \\
\le& Cs\la\va(T) \Vert y(\cdot,T)\Vert^2_{H^1(\OOO)} e^{2s\va(T)}
+ Cs\la\Vert y(\cdot,0)\Vert^2_{H^1(\OOO)}e^{2s}
\end{align*}
for all large $s>0$.

Now the same argument as in the proof of Theorem 1 can complete the proof
of Theorem 3.
$\blacksquare$
\section{Proof of Theorem 2}

The proof is based on the stability for $t_0>0$ and 
$$
u(x,0) = \int^0_T \ppp_tu(x,t) dt + u(x,T), \quad x\in \OOO.
$$

We apply the same arguments towards (3.8) for $\ppp_tu$ and 
$\ppp_t^2u$.
We set
$$
M_1:= \sum_{k=0}^2 \Vert \ppp_t^ku(\cdot,0)\Vert_{H^1(\OOO)},
\quad D:= \sum_{k=0}^2 \Vert \ppp_t^ku(\cdot,T)\Vert_{H^1(\OOO)}.
$$
Then, arguing similarly to the proof of Theorem 1 to 
$\ppp_t^2u$, we obtain (3.8) for $\ppp_t^2u$ replacing $u$:
$$
\int_{\OOO\times (t_0, T)}
\left( \frac{1}{s}\vert \ppp_t^3u\vert^2
+ s\vert \ppp_t^2u\vert^2\right) dxdt
\le CsD^2 e^{2s\va(T)} + CsM_1^2e^{-2s\mu(t_0)},
$$
that is,    
$$
\Vert \ppp_t^2u\Vert^2_{L^2(\OOO\times (t_0,T))}
\le CD^2 e^{2s\va(T)} + CM_1^2e^{-2s\mu(t_0)}    \eqno{(4.1)}
$$
for all $s\ge s_0$: some constant and all $0<t_0<T$.

Since
$$
\ppp_tu(x,t_0) = \int^{t_0}_T \ppp_t^2u(x,t) dt + \ppp_tu(x,T), 
\quad x\in \OOO,
$$
we have
$$
\Vert \ppp_tu(\cdot,t_0)\Vert^2_{L^2(\OOO)}
\le C\Vert \ppp_t^2u\Vert^2_{L^2(\OOO\times (t_0,T))}
+ C\Vert \ppp_tu(\cdot,T)\Vert^2_{L^2(\OOO)}.      \eqno{(4.2)}
$$
We substitute (4.1) to (4.2), so that 
$$
\Vert \ppp_tu(\cdot,t_0)\Vert^2_{L^2(\OOO)}
\le CD^2 e^{2s\va(T)} + CM_1^2 e^{-2s\mu(t_0)}
+ C\Vert \ppp_tu(\cdot,T)\Vert^2_{L^2(\OOO)}
$$
$$
\le C_1D^2 e^{2s\va(T)} + C_1M_1^2 e^{-2s\mu(t_0)}    \eqno{(4.3)}
$$
for all $s\ge s_0$: some constant and all $0<t_0<T$.
Here we used $\Vert \ppp_tu(\cdot,T)\Vert_{L^2(\OOO)}^2
\le D^2 \le CD^2e^{2s\va(T)}$.

Using 
$$
u(x,0) = \int^0_T \ppp_tu(x,t_0) dt_0 + u(x,T), \quad x\in \OOO,
$$
we obtain
\begin{align*}
& \int_{\OOO} \vert u(x,0)\vert^2 dx
\le 2\int_{\OOO}\left\vert \int^T_0 \ppp_tu(x,t_0) dt_0\right\vert^2
dx + 2\int_{\OOO} \vert u(x,T)\vert^2 dx\\
\le& C\int^T_0 \Vert \ppp_tu(\cdot,t_0)\Vert^2_{L^2(\OOO)} dt_0
+ C\Vert u(\cdot,T)\Vert^2_{L^2(\OOO)}.
\end{align*}
Substitution of (4.3) yields 
$$
\Vert u(\cdot,0)\Vert^2_{L^2(\OOO)}
\le C_1D^2\int^T_0 e^{2s\va(T)} dt_0
+ C_1M_1^2\int^T_0 e^{-2s\mu(t_0)} dt_0
+ C\Vert u(\cdot,T)\Vert^2_{L^2(\OOO)},    
$$
that is,
$$
\Vert u(\cdot,0)\Vert^2_{L^2(\OOO)}
\le C_1D^2\, e^{C_2s} + C_1M_1^2\int^T_0 e^{-2s\mu(t_0)} dt_0.
                                          \eqno{(4.4)}
$$
We calculate $\int^T_0 e^{-2s\mu(t_0)} dt_0$ as follows.
Changing the variables $t_0 \mapsto \xi$ by 
$\xi:= \mu(t_0) = e^{\la t_0} - 1$ and setting $\www{T}:= e^{\la T}-1$, 
we have
$$
\int^T_0 e^{-2s\mu(t_0)} dt_0 = \frac{1}{\la}\int^{\www{T}}_0
e^{-2s\xi}\frac{1}{1+\xi} d\xi
\le \frac{1}{\la}\left[ \frac{e^{-2s\xi}}{2s} \right]^{\xi=0}
_{\xi=\www{T}}
\le \frac{1}{2\la s}.
$$

Hence, (4.4) yields
$$
\Vert u(\cdot,0)\Vert^2_{L^2(\OOO)}
\le C_3D^2\, e^{C_2s} + \frac{C_3}{s}M_1^2                    \eqno{(4.5)}
$$
for all $s \ge s_0$.

We choose $s>0$ suitably for making the right-hand side small.
Taking into consideration that the first term $e^{C_2s}$ increases and the 
second decreases with $\frac{1}{s}$, we may choose $s>0$ like order
$\left(\log \frac{1}{D}\right)^{-1}$.  Precisely, we argue as follows.

First setting $C_4:= C_3e^{C_2s_0}$ and replacing $s:= s+s_0$. we can 
obtain (4.5) for all $s>0$.  Without loss of generality, we can assume
$D < 1$.  Setting
$$
s = \left( \log \frac{1}{D} \right)^{\alpha} > 0
$$
with $0<\alpha<1$, we have
\begin{align*}
& e^{C_2s}D^2 
= \exp\left( C_2\left( \log \frac{1}{D} \right)^{\alpha} \right)D^2\\
=& \exp\left( C_2\left( \log \frac{1}{D} \right)^{\alpha} \right)
e^{2\log D}
= \exp\left( -2\left( \log \frac{1}{D} \right)
+ C_2\left( \log \frac{1}{D} \right)^{\alpha} \right).
\end{align*}
Since we can find a constant $C_5>0$ such that  
$$
e^{-2\eta + C_2\eta^{\alpha}} \le \frac{C_5}{\eta^{\alpha}}\quad 
\mbox{for all $\eta>0$}
$$
by $\alpha < 1$, we see
$$
e^{C_2s}D^2 \le C_5\left( \log \frac{1}{D} \right)^{-\alpha},
$$
and 
$$
C_4D^2e^{C_2s} + \frac{C_4}{s}M_1^2
\le C_4C_5\left( \log \frac{1}{D} \right)^{-\alpha} 
+ C_3M_1^2\left( \log \frac{1}{D} \right)^{-\alpha} .
$$
Thus the proof of Theorem 2 is complete.
$\blacksquare$
%
%
%
%\section{Concluding remarks}

{\bf Acknowledgements.}
The work was supported by Grant-in-Aid for Scientific Research (A) 20H00117 
of Japan Society for the Promotion of Science.

\end{document}